\documentclass[a4paper,english,11pt]{article}

\usepackage[english]{babel}
\usepackage{theorem,epsfig,subfigure,array,mathe}
\usepackage{varioref,float,amssymb,amsmath,float}
\usepackage[latin1]{inputenc}
\newtheorem{theor}{Theorem}
\newtheorem{lemma}[theor]{Lemma}

\newtheorem{remrk}[theor]{Remark}

\newcommand{\me}{1/2}
\newcommand{\D}{\protect\displaystyle}

\newcommand{\ipl}{\langle}
\newcommand{\ipr}{\rangle}
\newcommand{\eop}{\hspace*{\fill} $\blacksquare$}
\begin{document}
\pagestyle{headings}
\markboth{Iterative methods for ill-posed problems modeled by PDE's}
{J.Baumeister and A.Leit\~ao}
% -------------------------------------------------------------------------- %
%                                  Title Page                                %
% -------------------------------------------------------------------------- %
\title{On iterative methods for solving ill-posed problems 
modeled by partial differential equations}

\setcounter{footnote}{1}
\author{J.Baumeister\footnote{Fachbereich Mathematik, Johann Wolfgang Goethe 
Universit\"at, Robert--Mayer--Str. 6--10, 60054 Frankfurt am Main, Germany 
(baumeist@math.uni-frankfurt.de)}
\ and 
A.Leit\~ao\footnote{Department of Mathematics, Federal University 
of Santa Catarina, P.O. Box 476, 88010-970 Florian\'opolis, Brazil 
(aleitao@mtm.ufsc.br)} }

\maketitle

\begin{abstract}
        We investigate the iterative methods proposed by Maz'ya and Kozlov 
(see [KM1], [KM2]) for solving ill-posed inverse problems modeled by partial 
differential equations. We consider linear evolutionary problems of 
elliptic, hyperbolic and parabolic types. Each iteration of the analyzed 
methods consists in the solution of a well posed problem (boundary value 
problem or initial value problem respectively). The iterations are described 
as powers of affine operators, as in [KM2]. We give alternative convergence 
proofs for the algorithms by using spectral theory and the fact that the 
linear parts of these affine operators are non-expansive with additional 
functional analytical properties (see [Le1,2]). Also problems with noisy 
data are considered and estimates for the convergence rate are obtained 
under {\em a priori} regularity assumptions on the problem data.
\end{abstract}
%
%
%
% -------------------------------------------------------------------------- %
%                                  Section 1                                 %
% -------------------------------------------------------------------------- %
\setcounter{footnote}{1}
\section{Introduction}
%
%
%
%----------------------------------------------------------------------------%
\subsection{Main results}

        We present new convergence proofs for the iterative algorithms 
proposed in [KM2] using a functional analytical approach, were each 
iteration is described using powers of an affine operator $T$. The key of 
the proof is to choose a correct topology for the Hilbert space were the 
iteration takes place and to prove that $T_l$, the linear component of $T$, 
is a {\em regular asymptotic}, {\em non-expansive} operator (other 
properties of $T_l$ such as positiveness, self-adjointness and injectivity 
are also verified). The converse is also proved, i.e. if an iterative 
procedure converges, the limit point is the solution of the respective 
problem.

          The convergence rate of the iterative method can be estimated when 
we make appropriate regularity assumptions on the problem data. In the last 
section some numerical experiments are presented, were we test the algorithm 
performance for linear elliptic, hyperbolic and parabolic ill-posed problems.

        The iterative procedures discussed in this paper were presented in 
[KM2] and also treated via semi groups in [Bas]. The iterative procedure for 
elliptic Cauchy problems defined in domains of more general type is 
discussed in [KM1], [Le1,2] and [JoNa]. The iterative procedure concerning 
parabolic problems is also treated in [Va].
%
%
%
%----------------------------------------------------------------------------%
\subsection{Preliminaries} \label{ssec-prelim}

%========================================
\subsubsection{On non-expansive operators} \label{sssec-prelim-fanalysis}

        Let $H$ be a separable Hilbert space endowed with an inner product 
$\ipl \cdot,\cdot \ipr$ and norm $\| \cdot \|$. A linear operator $T: H \to 
H$ is called {\em non-expansive} if $\| T \| \leq 1$.

        An operator $T: H \to H$ is said to be {\em regular asymptotic} in 
$x \in H$ if
\[ \lim_{k\to\infty} \| T^{k+1}(x) - T^k(x) \| = 0  \]
holds true. If the above property holds for every $x \in H$, we say that 
$T$ is {\em regular asymptotic} in $H$.

        Next we formulate the results used to prove the convergence of 
the iterative algorithms analyzed in this paper.

\begin{lemma} \label{satz-nexp-asreg}

        Let $T:H \to H$ be a linear non-expansive operator. With $\Pi$ we 
denote the orthogonal projector defined on $H$ onto the null space of 
$(I-T)$. The following assertions are equivalent: \\[1ex]
a) $T$ is regular asymptotic in $H$; \\[1ex]
b) $\D\lim_{k\to\infty} T^k x = \Pi x$ for all $x \in H$.
\end{lemma}

        A proof of this lemma (even in a more general framework) can be found 
in [Je] (see also [Le1] and the references cited therein).

\begin{lemma} \label{satz-konverg-nexp-asreg}

        Let $T:H \to H$ be a linear, non-expansive, regular asymptotic 
operator such that 1 is not an eigenvalue of $T$.%
\footnote{The set of all eigenvalues of a linear operator $T$ is denoted by 
$\sigma_p(T)$.}
Given $z \in H$ define $S: H \ni x \mapsto Tx + z \in H$. Then for every 
$x_0 \in H$ the sequence $\{S^k x_0\}$ converges to the uniquely determined 
solution of the fixed point equation $Sx = x$. \\
Proof:
Let $\bar{x} \in H$ be the solution of $S\bar{x} = \bar{x}$. Defining $x_k := 
S^k x_0$ and $\eps_k := \bar{x} - x_k$ one can easily see that $\eps_{k+1} = 
T \eps_k$, $k \in \N$. Lemma~\ref{satz-nexp-asreg} allow us to conclude that 
$\lim_k \eps_k = \Pi \eps_0$. From the hypothesis we have $\ker(I-T) = \{0\}$, 
and the lemma follows. \eop
\end{lemma}

        In the next lemma we present a sufficient condition for an operator 
to be non-expansive and regular asymptotic. For convenience of the reader we 
include here the proof (see [KM2]).

\begin{lemma} \label{satz-mazya-bedingung}

        Let $T$ be a bounded linear operator in $H$ such that for $c>0$
\begin{equation} \label{gl-mazya-bedingung}
 \| (I-T)x \|^2 \ \leq \ c( \| x \|^2 - \| Tx \|^2 )\, ,\ 
   \forall x \in H ,
\end{equation}
holds true. Then $T$ is non-expansive and regular asymptotic in $H$. \\
Proof:
The non-expansivity of $T$ follows directly from the inequality
\[ 0 \ \leq \ c^{-1} \| (I-T)x \|^2 \ \leq \ \| x \|^2 - \| Tx \|^2 ,\ 
   \forall x \in H . \]
Now take $x_0 \in H$. Since $\|T\| \leq 1$, the sequence $\| T^k x_0 \|^2$ 
is non-increasing, from what we conclude that $\lim_k( \| T^k x_0 \|^2 - 
\| T^{k+1} x_0 \|^2 ) = 0$. Note that from (\ref{gl-mazya-bedingung}) follows
\[ \| T^k x_0 - T^{k+1} x_0 \|^2 \ \leq \ c( \| T^k x_0 \|^2 - 
   \| T^{k+1} x_0 \|^2 ) . \]
Putting all together one can see that $T$ is regular asymptotic in $H$. \eop
\end{lemma}

        Equivalent to the condition (\ref{gl-mazya-bedingung}) in 
Lemma~\ref{satz-mazya-bedingung} is the following one%
\footnote{Clearly, condition (\ref{gl-mazya-bedingung}) is only sufficient 
for $T$ being non-expansive and regular asymptotic in $H$.}
\begin{equation} \label{equiv-mazya-bedingung}
\ipl (I-T) x , x \ipr \ \geq \ \frac{c+1}{2c} \|(I-T) x\|^2 ,\ 
\forall x\in H ,
\end{equation}
as the following line suggests (see [KM2])
\[  \| x \|^2 \ = \ \|Tx\|^2 - \|(I-T) x\|^2 + 2 \ipl x , (I-T)x \ipr
\, ,\ \forall x \in H . \]

%========================================
\subsubsection{On function spaces} \label{sssec-prelim-fspaces}

        Let $\Omega \subset \R^n$ be an open, bounded set with smooth 
boundary and let $A$ be a positive, self-adjoint, unbounded operator densely 
defined on the Hilbert space $H := L_2(\Omega)$. Let $E_\lbd$, $\lbd \in \R$, 
denote the resolution of the identity associated to $A$, i.e.
\[ \ipl A \vphi, \psi \ipr =
   \int_{\lbd\in\R} \lbd \, d \ipl E_\lbd\vphi,\psi \ipr =
   \int_0^\infty \lbd \, d \ipl E_\lbd\vphi,\psi \ipr , \]
for $\vphi \in D(A)$, the domain of $A$, and $\psi \in H$. Note that 
given $f \in C(\R^+)$ we can define the operator $f(A)$ on $H$ by setting
\[ \ipl f(A) \vphi, \psi \ipr :=
   \int_0^\infty f(\lbd) \, d \ipl E_\lbd \vphi,\psi \ipr , \]
for every $\vphi \in D(f(A))$ and $\psi \in H$, where the domain of $f(A)$ is 
defined by
\[ D(f(A)) := \{ \vphi \in H \, | \,
   \int_0^\infty f(\lbd)^2 d \ipl E_\lbd \vphi,\vphi \ipr < \infty \} . \]
Now we are ready to construct a family of Hilbert spaces ${\cal H}^s(\Omega)$, 
$s \geq 0$, as the domain of definition of the powers of $A$%
\footnote{For simplicity we may write ${\cal H}^s$ instead of 
${\cal H}^s(\Omega)$.}
\begin{equation} \label{def-Hs-raeume}
{\cal H}^s(\Omega)  :=  \{ \vphi \in H \, | \, \|\vphi\|_s := \left(
 \int_0^\infty (1+\lbd^2)^s d \ipl E_\lbd \vphi,\vphi \ipr \right)^{\me}
 < \infty \} .
\end{equation}
The Hilbert spaces ${\cal H}^{-s}(\Omega)$ (with $s > 0$) are defined by 
duality:%
\footnote{Alternatively one can define ${\cal H}^{-s}(\Omega)$ as the 
completion of $H$ in the (-s)--norm defined in (\ref{def-Hs-raeume}).}
${\cal H}^{-s} := ({\cal H}^{s})'$. It follows from the definition that 
${\cal H}^0(\Omega) = H$. It can also be proved that the embedding 
${\cal H}^r(\Omega) \hookrightarrow {\cal H}^s(\Omega)$ is dense and compact 
for $r>s$ (see [LiMa] Chapter 1).

        An interesting example is $A = (-\Delta)^{\me}$, where $\Delta$ is 
the Laplace--Beltrami operator on $\Omega$. In this particular case we have 
the identity ${\cal H}^s(\Omega) = H^{2s}_0(\Omega)$, where $H^{s}_0(\Omega)$ 
is the Sobolev space of index $s$ according to the definition of Lions and 
Magenes (see [LiMa] pp. 54). One should note that functions in ${\cal H}^s 
(\Omega)$ satisfy null boundary conditions in the sense of the trace operator.

        Given $T>0$ we define the spaces $L_2(0,T; {\cal H}^s(\Omega))$ of 
functions $u: (0,T) \ni t \mapsto u(t) \in {\cal H}^s(\Omega)$. These are 
normed spaces if we consider
\[ \|u\|_{2;0,T;s} \ := \ 
   \left( \int_0^T \|u(t)\|_s^2 \, dt \right)^{\me} , \]
as a norm in $L_2(0,T; {\cal H}^s(\Omega))$. Finally, we define the spaces 
$C(0,T; {\cal H}^s(\Omega))$ of continuous functions $u: [0,T] \ni t 
\mapsto u(t) \in {\cal H}^s(\Omega)$. The norm on these spaces is given by
\[ \|u\|_{\infty;0,T;s} \ := \ \sup_{t\in[0,T]} \|u(t)\|_s . \]
%
%
%
% -------------------------------------------------------------------------- %
%                                  Section 2                                 %
% -------------------------------------------------------------------------- %
\setcounter{footnote}{1}
\section{The ill-posed problems} \label{sec-ill-posed-probl}

        Let the operator $A$ with discrete spectrum, the set $\Omega$ and 
the Hilbert spaces ${\cal H}^s(\Omega)$ be defined as in 
section~\ref{ssec-prelim}. In the next three paragraphs we formulate the 
ill-posed problems that are discussed in this article.
%
%
%
%----------------------------------------------------------------------------%
\subsection{The elliptic problem:} \label{ssec-elliptic-probl-def}

        Given functions $(f,g) \in {\cal H}^{\me}(\Omega) \times 
{\cal H}^{-\me}(\Omega)$, find $u \in (V_e, \|\cdot\|_{V_e})$, where
\begin{eqnarray*}
V_e := L_2(0,T; {\cal H}^1(\Omega)) \\
\|u\|_{V_e} := \left(
      \int_0^T ( \|u(t)\|_1^2 + \|\partial_t u(t)\|_0^2 )\, dt \right)^{\me} ,
\end{eqnarray*}
that satisfies \\[2ex]
$ (P_e) \hskip2cm \left\{  \begin{array}{l}
            (\partial_t^2 - A^2) u = 0,\ {\rm in}\ (0,T) \times \Omega \\
            u(0,x)                 = f(x),\
            \partial_t u(0,x)      = g(x),\ x \in \Omega .
         \end{array} \right. $ \\[2ex]
Note that if $u \in V_e$, then $\partial_t u \in L_2(0,T; H)$ and appropriate 
trace theorems (see [LiMa]) guarantee that $u(0), u(T) \in 
{\cal H}^{\me}(\Omega)$ and $\partial_t u(0), \partial_t u(T) \in 
{\cal H}^{-\me}(\Omega)$.

        In this problem we are mostly interested in the value of $u$ for 
$t=T$, i.e. $u(T,x)$ and $\partial_t u(T,x)$, $x\in\Omega$. This elliptic 
initial value problem (also called Cauchy problem) is not well posed in the 
sense of Hadamard (see [Bau]). This follows from the general representation 
of the solution of $(P_e)$ given by
\begin{equation}
 u(t) = \cosh(At) f + \sinh(At) A^{-1} g .
\end{equation}
One can construct a sequence of Cauchy data $(f_k,g_k) = (0,g_k)$ using the 
eigenfunctions of $A$, such that $(f_k,g_k)$ converge to zero in 
${\cal H}^{\me} \times {\cal H}^{-\me}$ while the norm of the solutions 
$\|u_k\|_{V_e}$ do not.
%
%
%
%----------------------------------------------------------------------------%
\subsection{The hyperbolic problem:} \label{ssec-hyperbolic-probl-def}

        Given functions $f, g \in {\cal H}^1(\Omega)$, find $u \in (V_h, 
\|\cdot\|_{V_h})$, where
\begin{eqnarray*}
V_h := \{ v \in C(0,T; {\cal H}^1(\Omega)) \, | \,
         \partial_t u \in C(0,T; H) \} \\
\|u\|_{V_h} := \sup_{t\in[0,T]} \left(
                        \|u(t)\|_1^2 + \|\partial_t u(t)\|_0^2 \right)^{\me} ,
\end{eqnarray*}
that satisfies \\[2ex]
$ (P_h) \hskip2cm \left\{  \begin{array}{l}
             (\partial_t^2 + A^2) u = 0,\ {\rm in}\ (0,T) \times \Omega \\
             u(0,x)                 = f(x),\
             u(T,x)                 = g(x),\ x \in \Omega .
           \end{array} \right. $ \\[2ex]
Note that if $u \in V_h$, then $u(0), u(T) \in {\cal H}^1(\Omega)$ and 
$\partial_t u(0), \partial_t u(T) \in H$.

        Let's assume that the numbers $k\pi/T$, $k=1,2,\dots$ are not 
eigenvalues of $A$.%
\footnote{If this condition is not satisfied, one can easily see that problem 
$(P_h)$ is not uniquely solvable.}
Then this hyperbolic (Dirichlet) boundary value problem is ill-posed if the 
distance from the set $M := \{ k\pi/T ; k \in \N \}$ to $\sigma(A)$ -- the 
spectrum of $A$ -- is zero. To see this, we take $\lbd_k \in \sigma(A)$ with 
$\lim_k$ dist$(\lbd_k,M) = 0$ and $g_k$ the respective (normalized) 
eigenfunctions. Solving problem $(P_h)$ for the data $(f,g)= (0,g_k)$ one 
obtains respectively the solutions
\begin{equation}
 u_k(t) = \sin(At)\sin(AT)^{-1} g_k
        = \sin(\lbd_k t) \sin(\lbd_k T)^{-1} g_k ,
\end{equation}
which happens to be unbounded in $V_h$.
%
%
%
%----------------------------------------------------------------------------%
\subsection{The parabolic problem:} \label{ssec-parabolic-probl-def}

        Given a function $f \in H = L_2(\Omega)$ find $u \in (V_p, 
\|\cdot\|_{V_p})$, where
\begin{eqnarray*}
V_p := L_2(0,T; {\cal H}^1(\Omega)) \\
\|u\|_{V_p} := \left(
   \int_0^T ( \|u(t)\|_1^2 + \|\partial_t u(t)\|_{-1}^2 )\, dt \right)^{\me} ,
\end{eqnarray*}
that satisfies \\[2ex]
$ (P_p) \hskip2cm \left\{ \begin{array}{l}
             (\partial_t + A^2) u = 0,\ {\rm in}\ (0,T)\times\Omega \\
             u(T,x)               = f(x), \ x \in \Omega
          \end{array} \right. $ \\[2ex]
Note that if $u \in V_p$, then $u(0), u(T) \in H$.

         Problem $(P_p)$ corresponds to the well known problem of solving the 
heat equation backwards in time, which is known to be (severely) ill-posed. 
This follows from the general representation of the solution of $(P_p)$ 
given by
\begin{equation} \label{loesung-parab-probl}
 u(t) = \exp(A^2(T-t)) f .
\end{equation}
Again using the eigenfunctions of $A$, one can construct a sequence of data 
$f_k$ converging to zero in $H$ while the norm of the solutions 
$\|u_k\|_{V_p}$ do not.
%
%
%
% -------------------------------------------------------------------------- %
%                                  Section 3                                 %
% -------------------------------------------------------------------------- %
\setcounter{footnote}{1}
\section{Description of the methods} \label{sec-methoden}
%
%
%
%----------------------------------------------------------------------------%
\subsection{The iterative procedure for the elliptic problem}
\label{ssec-it-proc-elliptic}

        Consider problem $(P_e)$ with data $(f,g) \in {\cal H}^{\me}(\Omega) 
\times {\cal H}^{-\me}(\Omega)$. Given any initial guess $\vphi_0 \in 
{\cal H}^{-\me}(\Omega)$ for $\partial_t u(T)$ we improve it by solving 
the following mixed boundary value problems (BVP) of elliptic type:
$$ \left\{  \begin{array}{l}
     (\partial_t^2 - A^2) v = 0,\ {\rm in}\ (0,T) \times \Omega \\
     v(0) = f,\ \ \partial_t v(T) = \vphi_0
   \end{array} \right. \hskip1cm
   \left\{  \begin{array}{l}
     (\partial_t^2 - A^2) w = 0,\ {\rm in}\ (0,T) \times \Omega \\
     \partial_t w(0) = g,\ \ w(T) = v(T)
     
   \end{array} \right. $$
and defining $\vphi_1 := \partial_t w(T)$. Each one of the mixed BVP's 
above has a solution in $V_e$ and consequently $\vphi_1 \in {\cal H}^{-\me} 
(\Omega)$. Setting $\vphi_0 := \vphi_1$ and repeating this procedure we 
construct a sequence $\{\vphi_k\}$ in ${\cal H}^{-\me}(\Omega)$.

        Our assumptions on the operator $A$ allow the determination of the 
exact solutions $v$ and $w$ of the above problems, which are given by
\[ v(t) = \sinh(At) \cosh(AT)^{-1} \, A^{-1} \vphi_0 + 
                \cosh(A(t-T)) \cosh(AT)^{-1} f , \]
\[ w(t) = \cosh(At) \cosh(AT)^{-1} v(T) + 
                \sinh(A(t-T)) \cosh(AT)^{-1} A^{-1} g . \]
Finally, we can write
\[ \vphi_1 = \partial_t w(T) = \tanh(AT)^2 \vphi_0 +
                               \sinh(AT) \cosh(AT)^{-2} A f +
                               \cosh(AT)^{-1} g . \]
Now, defining the affine operator $T_e: {\cal H}^{-\me}(\Omega) \to 
{\cal H}^{-\me}(\Omega)$ by
\begin{equation} \label{Te-oper-def}
 T_e(\vphi) := \tanh(AT)^2 \vphi + z_{f,g} ,
\end{equation}
with $z_{f,g} := \sinh(AT) \cosh(AT)^{-2} A f + \cosh(AT)^{-1} g$, the 
iterative algorithm can be rewritten as
\begin{equation} \label{elliptic-iter-def}
\vphi_k = T_e(\vphi_{k-1}) = T_e^k(\vphi_0) =
\tanh(AT)^{2k} \vphi_0 + \sum_{j=0}^{k-1} \tanh(AT)^{2j} z_{f,g} .
\end{equation}
%
%
%
%----------------------------------------------------------------------------%
\subsection{The iterative procedure for the hyperbolic problem}
\label{ssec-it-proc-hyperbolic}

        Let's now consider problem $(P_h)$ with data $f, g \in {\cal H}^1 
(\Omega)$. Given any initial guess $\vphi_0 \in H$ for $\partial_t u(0)$ we 
improve it by solving the following initial value problems (IVP) of 
hyperbolic type:%
\footnote{The second problem is considered with reversed time.}
$$ \left\{  \begin{array}{l}
     (\partial_t^2 + A^2) v = 0,\ {\rm in}\ (0,T) \times \Omega \\
     v(0) = f,\ \ \partial_t v(0) = \vphi_0
   \end{array} \right. \hskip1cm
   \left\{  \begin{array}{l}
     (\partial_t^2 + A^2) w = 0,\ {\rm in}\ (0,T) \times \Omega \\
     w(T) = g,\ \ \partial_t w(T) = \partial_t v(T)
   \end{array} \right. $$
and defining $\vphi_1 := \partial_t w(0)$. Each one of the mixed IVP's 
above has a solution in $V_h$ and consequently $\vphi_1 \in H$. Repeating 
this procedure we construct a sequence $\{\vphi_k\}$ in $H$.

        As in section~\ref{ssec-it-proc-elliptic}, the assumptions on the 
operator $A$ allow the determination of the exact solutions $v$ and $w$ 
of the above problems. In fact we have
\[ v(t) = \cos(At) f + \sin(At) A^{-1} \vphi_0 , \]
\[ w(t) = \cos(A(t-T)) \, g + \sin(A(t-T)) A^{-1} \partial_t v(T) . \]
Finally, we can write
\[ \vphi_1 = \partial_t w(0) = \cos(AT)^2 \vphi_0 -
               \cos(AT)\sin(AT) A f \ + \ \sin(AT) g \]
and defining the affine operator $T_h: H \to H$ by
\begin{equation} \label{Th-oper-def}
 T_h(\vphi) := \cos(AT)^2 \vphi + z_{f,g} ,
\end{equation}
with $z_{f,g} := -\cos(AT)\sin(AT) A f + \sin(AT) g$, the iterative 
method can be rewritten as
\begin{equation} \label{hyperbolic-iter-def}
\vphi_k = T_h(\vphi_{k-1}) = T_h^k(\vphi_0) =
\cos(AT)^{2k} \vphi_0 + \sum_{j=0}^{k-1} \cos(AT)^{2j}\, z_{f,g} .
\end{equation}
%
%
%
%----------------------------------------------------------------------------%
\subsection{The iterative procedure for the parabolic problem}
\label{ssec-it-proc-parabolic}

        We consider problem $(P_p)$ with data $f \in H$. Define $\bar{\lbd} 
:= \inf\{ \lbd; \lbd \in \sigma(A) \}$ and chose a positive parameter 
$\gamma$ such that $\gamma < 2 \exp(\bar{\lbd}^2 T)$. Now, given $\vphi_0 \in 
H$ an initial guess for $u(0)$, the method consists in first solving the IVP 
of parabolic type:
\[ \left\{  \begin{array}{l}
     (\partial_t + A^2) v_0 = 0,\ {\rm in}\ (0,T) \times \Omega \\
     v_0(0) = \vphi_0
   \end{array} \right. \]
Then we solve for $k \geq 1$ the sequence of IVP's:
\[ \left\{  \begin{array}{l}
     (\partial_t + A^2) v_k = 0,\ {\rm in}\ (0,T) \times \Omega \\
     v_k(0) = v_{k-1}(0) - \gamma (v_{k-1}(T) - f)
   \end{array} \right. \]
The sequence $\{\vphi_k\}$ is defined by $\vphi_k := v_k(0) \in H$. 
Note that the analytic solutions of the above problems are given by
\[  v_k(t) = \exp(-A^2 t) \vphi_k , \]
and we obtain
\[ \vphi_{k+1} = (I - \gamma \exp(-A^2 T)) \vphi_k + \gamma f . \]
Now, we define the affine operator $T_p: H \to H$ by
\begin{equation} \label{Tp-oper-def}
T_p(\vphi) := (I - \gamma \exp(-A^2 T)) \vphi + z_f ,
\end{equation}
with $z_f :=  \gamma f$, and we are able to rewrite the iterative algorithm as
\begin{eqnarray}
\vphi_k & = & T_p(\vphi_{k-1}) = T_p^k(\vphi_0) \nonumber \\
        & = & (I - \gamma \exp(-A^2 T))^k \vphi_0 + \sum_{j=0}^{k-1}
              (I - \gamma \exp(-A^2 T))^j z_f . \label{parabolic-iter-def}
\end{eqnarray}
%
%
%
% -------------------------------------------------------------------------- %
%                                  Section 4                                 %
% -------------------------------------------------------------------------- %
\setcounter{footnote}{1}
\section{Analysis of the methods}
%
%
%
%----------------------------------------------------------------------------%
\subsection{The elliptic case} \label{ssec-ellip-analyse}

        The linear part of the affine operator $T_e$ defined in 
(\ref{Te-oper-def}) is given by $T_{l,e} := \tanh(AT)^2$. We begin the 
discussion analyzing an important property of problem $(P_e)$.

\begin{lemma} \label{lemma-cauchy-kowal}
        Given $(f,g) \in {\cal H}^{\me} \times {\cal H}^{-\me}$, problem 
$(P_e)$ has at most one solution in $V_e$. \\
Proof:
This result is a generalization of the Cauchy--Kowalewsky theorem. A complete 
proof can be found in [Le1]. \eop
\end{lemma}

        From Lemma~\ref{lemma-cauchy-kowal} follows that if problem $(P_e)$ 
has a solution $u \in V_e$, then it's Neumann trace $\bar{\vphi} := 
\partial_t u(T)$ solves the equation $T_e \bar{\vphi} = \bar{\vphi}$. The 
objective of the iterative method in section~\ref{ssec-it-proc-elliptic} is to 
find a solution of this fixed point equation. The ill-posedness of problem 
$(P_e)$ can be recognized in the fact that 1 belongs to continuous spectrum 
of $T_{l,e}$, as one can see in the next lemma.

\begin{lemma} \label{lemma-Tle-eigenschaft}
        The linear operator $T_{l,e}: {\cal H}^{-\me} \to {\cal H}^{-\me}$ 
is positive, self-adjoint, injective, non-expansive, regular asymptotic and 
1 is not an eigenvalue of $T_{l,e}$.
Further $T_{l,e}$ satisfies the condition (\ref{gl-mazya-bedingung}). \\
Proof:
The injectivity follows promptly from Lemma~\ref{lemma-cauchy-kowal}. The 
properties: positiveness, self-adjointness and $1 \not\in \sigma_p(T_{l,e})$ 
follow from the definition of $T_{l,e}$ together with the assumptions on $A$ 
made in section~\ref{sssec-prelim-fspaces} (remember we required in 
section~\ref{sec-ill-posed-probl} that $\sigma(A)$ is discrete).

        In order to prove that $T_{l,e}$ is non-expansive and regular 
asymptotic, it is enough to verify the condition (\ref{gl-mazya-bedingung}) 
(see Lemma~\ref{satz-mazya-bedingung}). It's easy to see that $T_{l,e}$ 
satisfies this condition with $c=1$, if $\sigma(T_{l,e}) \in [0,1]$. One 
should note that this last property was already proved above. \eop
\end{lemma}

        In the next theorem we discuss the convergence of the algorithm 
described in section~\ref{ssec-it-proc-elliptic}.

\begin{theor} \label{satz-converg-ellip-verfahren}
        Let $T_e$ be the operator defined in (\ref{Te-oper-def}) and 
$T_{l,e}$ it's linear part. If problem $(P_e)$ in 
section~\ref{ssec-elliptic-probl-def} is consistent%
\footnote{This means that it has a corresponding solution $u \in V_e$.}
for the data $(f,g)$, then the sequence $\{\vphi_k\}$ defined in 
(\ref{elliptic-iter-def}) converges to $\partial_t u(T)$ in the 
norm of ${\cal H}^{-\me}(\Omega)$. \\
Proof:
Follows from Lemma~\ref{lemma-Tle-eigenschaft} and 
Lemma~\ref{satz-konverg-nexp-asreg} with $z := z_{f,g}$, $T := T_{l,e}$ 
and $S := T_e$. \eop
\end{theor}

        The converse of Theorem~\ref{satz-converg-ellip-verfahren} is also 
true, i.e. if the sequence $\{\vphi_k\}$ in (\ref{elliptic-iter-def}) 
converges in ${\cal H}^{-\me}(\Omega)$, it converges to the solution of 
$(P_e)$.

\begin{theor}
        If the sequence $\{\vphi_k\}$ defined in (\ref{elliptic-iter-def}) 
converges, say to $\bar{\vphi}$, then problem $(P_e)$ is consistent for the 
Cauchy data $(f,g)$ and it's solution $u \in V_e$ satisfies $\partial_t 
u(T) = \bar{\vphi}$. \\
Proof:
If $\lim_k \vphi_k = \bar{\vphi}$, then $T_e \bar{\vphi} = \bar{\vphi}$. 
Taking $\vphi_0=\bar{\vphi}$ in the mixed BVP's of 
section~\ref{ssec-it-proc-elliptic} we see that the functions $v$, $w$ satisfy 
the same boundary conditions (Dirichlet and Neumann conditions, respectively) 
at $t=T$. From Lemma~\ref{lemma-cauchy-kowal} we must have $v=w$ and one can 
see that $u := v = w$ is the solution of $(P_e)$, the identity 
$\partial_t u(T) = \bar{\vphi}$ being obvious. \eop
\end{theor}
%
%
%
%----------------------------------------------------------------------------%
\subsection{The hyperbolic case} \label{ssec-hyperb-analyse}

        The linear part of the affine operator $T_h$ defined in 
(\ref{Th-oper-def}) is given by $T_{l,h} := \cos(AT)^2$. We start the 
discussion proving some properties of this operator.

\begin{lemma} \label{lemma-Tlh-eigenschaft}
        The linear operator $T_{l,h}: H \to H$ is positive, self-adjoint, 
injective, non-expansive, regular asymptotic and 1 is not an eigenvalue 
of $T_{l,h}$. Further $T_{l,h}$ satisfies the condition 
(\ref{gl-mazya-bedingung}). \\
Proof:
The injectivity follows from the assumption $\{ k\pi/T; k\in\N \} \cap 
\sigma(A) = \emptyset$. The properties: positiveness, self-adjointness and 
$1 \not\in \sigma_p(T_{l,h})$ are proved like in 
Lemma~\ref{lemma-Tle-eigenschaft}.

        Again we use Lemma~\ref{satz-mazya-bedingung} to prove that $T_{l,h}$ 
is non-expansive and regular asymptotic. Since $\sigma(T_{l,h}) \in [0,1]$, 
the condition (\ref{gl-mazya-bedingung}) is obtained analogous as in 
Lemma~\ref{lemma-Tle-eigenschaft}. \eop
\end{lemma}

        From Lemma~\ref{lemma-Tlh-eigenschaft} follows that if problem $(P_h)$ 
has a solution $u \in V_h$, then it's Neumann trace $\bar{\vphi} := 
\partial_t u(0)$ solves the equation $T_h \bar{\vphi} = \bar{\vphi}$. Just 
like in the elliptic case (see section~\ref{ssec-ellip-analyse}) the 
objective of the method in section~\ref{ssec-it-proc-hyperbolic} is to 
approximate the solution of this fixed point equation. The ill-posedness of 
problem $(P_h)$ reflects in the fact that 1 belongs to continuous spectrum 
of $T_{l,h}$ (see Lemma~\ref{lemma-Tlh-eigenschaft}). In the next theorem we 
discuss the convergence of the algorithm described in 
section~\ref{ssec-it-proc-hyperbolic}.

\begin{theor} \label{satz-converg-hyperb-verfahren}
        Let $T_h$ be the operator defined in (\ref{Th-oper-def}) and 
$T_{l,h}$ it's linear part. If problem $(P_h)$ in 
section~\ref{ssec-elliptic-probl-def} is consistent for the data $(f,g)$, 
then the sequence $\{\vphi_k\}$ defined in (\ref{hyperbolic-iter-def}) 
converges to $\partial_t u(0)$ in the norm of $H$. \\
Proof:
Follows from Lemma~\ref{lemma-Tlh-eigenschaft} and 
Lemma~\ref{satz-konverg-nexp-asreg} with $z := z_{f,g}$, $T := T_{l,h}$ 
and $S := T_h$. \eop
\end{theor}

        The converse of Theorem~\ref{satz-converg-hyperb-verfahren} is also 
true, i.e. if the sequence $\{\vphi_k\}$ in (\ref{hyperbolic-iter-def}) 
converges in $H$, it converges to the solution of $(P_h)$.

\begin{theor}
        If the sequence $\{\vphi_k\}$ defined in (\ref{hyperbolic-iter-def}) 
converges, say to $\bar{\vphi}$, then problem $(P_h)$ is consistent for the 
Cauchy data $(f,g)$ and it's solution $u \in V_h$ satisfies $\partial_t 
u(0) = \bar{\vphi}$. \\
Proof:
If $\lim_k \vphi_k = \bar{\vphi}$, then $T_h \bar{\vphi} = \bar{\vphi}$. 
Taking $\vphi_0=\bar{\vphi}$ in the IVP's of 
section~\ref{ssec-it-proc-hyperbolic} we see that the functions $v$, $w$ 
satisfy the same Neumann boundary conditions at $t=0$ and $t=T$. From 
Lemma~\ref{lemma-Tlh-eigenschaft} we must have $v=w$ and one can 
see that $u := v = w$ is the solution of $(P_h)$, the identity 
$\partial_t u(0) = \bar{\vphi}$ being obvious. \eop
\end{theor}
%
%
%
%----------------------------------------------------------------------------%
\subsection{The parabolic case} \label{ssec-parab-analyse}

        The linear part of the affine operator $T_p$ defined in 
(\ref{Tp-oper-def}) is given by $T_{l,p} := I - \gamma \exp(-A^2 T)$. First, 
we analyze an important property of problem $(P_p)$.

\begin{lemma} \label{lemma-parab-exist-eindeut}
        Given $f \in H$, problem $(P_p)$ has exactly one solution in $V_p$. \\
Proof:
This result is suggested by the general representation of the solution given 
in (\ref{loesung-parab-probl}). A complete proof can be found in [LiMa], 
Chapter~3.
\mbox{} \eop
\end{lemma}

        Just like in the other cases, the iterative method in 
section~\ref{ssec-it-proc-parabolic} approximates the solution 
of the corresponding fixed point equation $T_p \bar{\vphi} = \bar{\vphi}$, 
which is uniquely solved by the Dirichlet trace $\bar{\vphi} = u(0)$ of the 
solution $u \in V_p$ of $(P_p)$ (see Lemma~\ref{lemma-parab-exist-eindeut}). 
Next, we discuss some properties of $T_{l,p}$.

\begin{lemma} \label{lemma-Tlp-eigenschaft}
        The linear operator $T_{l,p}: H \to H$ is self-adjoint, non-expansive, 
regular asymptotic and 1 is not an eigenvalue of $T_{l,p}$. Further, if 
$\gamma < 2\exp(\tilde{\lbd}^2T)$, where $\tilde{\lbd} := ( \bar{\lbd}^2 - 
T^{-1} \ln 2 )^{\me}$, then $T_{l,p}$ is injective and satisfies the 
condition (\ref{gl-mazya-bedingung}). \\
Proof:
The self-adjointness follows follows from the definition of $T_{l,p}$. Since 
the inequality $ 0 < \gamma \exp(-\lbd^2T) < 2 \exp([\bar{\lbd}^2 - \lbd^2]T) 
< 2$ holds for every $\lbd \in \sigma(A)$, we have $\sigma_p (T_{l,p}) \in 
(-1,1)$ and the non-expansivity follows. Note that the property $1 \not\in 
\sigma_p(T_{l,p})$ was also proved.

        To prove the asymptotic regularity we take $\vphi \in H$ and write
$T_{l,p}^{k+1} \vphi - T_{l,p}^k \vphi = T_{l,p}^k \psi$, where $\psi := 
(T_{l,p} - I) \vphi \in H$. Since $\sigma_p(T_{l,p}) \in (-1,1)$, it follows 
that $\lim_k T_{l,p}^k \psi = 0$, for all $\psi \in H$.

        Now, if $\gamma$ satisfies the extra assumption, a simple calculation 
shows that $\sigma_p(T_{l,p}) \in (0,1)$. The injectivity follows immediately 
and the condition (\ref{gl-mazya-bedingung}) is proved analogous as in 
Lemma~\ref{lemma-Tle-eigenschaft}. \eop
\end{lemma}

        In the next theorem we discuss the convergence of the algorithm 
described in section~\ref{ssec-it-proc-parabolic}.

\begin{theor} \label{satz-converg-parab-verfahren}
        Let $T_p$ be the operator defined in (\ref{Tp-oper-def}) and 
$T_{l,p}$ it's linear part. Given $f \in H$, let $u \in V_p$ be the 
uniquely determined solution of problem $(P_p)$. Then the sequence 
$\{\vphi_k\}$ defined in (\ref{parabolic-iter-def}) converges to $u(0)$ 
in the norm of $H$. \\
Proof:
Follows from Lemma~\ref{lemma-Tlp-eigenschaft} and 
Lemma~\ref{satz-konverg-nexp-asreg} with $z := z_f$, $T := T_{l,p}$ 
and $S := T_p$. \eop
\end{theor}
%
%
%
% -------------------------------------------------------------------------- %
%                                  Section 5                                 %
% -------------------------------------------------------------------------- %
\setcounter{footnote}{1}
\section{Regularization} \label{sec-regularization}

        In order to regularize the algorithms proposed in 
section~\ref{sec-methoden} we make the following assumptions on the 
formulation of the respective problems:
\begin{itemize}
\item[$(H_e)$] Given the Cauchy data $(f_\eps,g_\eps) \in {\cal H}^{\me} 
\times {\cal H}^{-\me}$, there exist consistent Cauchy data $(f,g) \in 
{\cal H}^{\me} \times {\cal H}^{-\me}$ such that $\|f-f_\eps\|_{\me} + 
\|g-g_\eps\|_{-\me} \leq \eps$, where $\eps >0$.
\item[$(H_h)$] Given the Dirichlet data $(f_\eps, g_\eps) \in {\cal H}^1 
\times {\cal H}^1$, there exist consistent Dirichlet data $(f,g) \in 
{\cal H}^1 \times {\cal H}^1$ such that $\|f-f_\eps\|_1 + \|g-g_\eps\|_1 
\leq \eps$, where $\eps >0$.
\item[$(H_p)$] The given data $f_\eps \in H$ is such that $\|f-f_\eps\|_H \leq 
\eps$, where $f \in H$ is the Dirichlet trace at $t=T$ of the exact solution 
of $(P_p)$ and $\eps >0$.
\end{itemize}
The assumptions on the data made in $(H_e)$ and $(H_h)$ may look very 
restrictive. One would prefer $f_\eps \in H = L_2(\Omega)$ in $(H_e)$ and 
$(f_\eps,g_\eps) \in H \times H$ in $(H_h)$, since these represent measured 
data. Nevertheless $(H_e)$ and $(H_h)$ are naturally satisfied if we make 
stronger assumptions on the regularity of the solutions of the corresponding 
ill-posed problems. This fact is explained in

\begin{lemma} \label{lemma-fehler-glaetung}
        Let $f \in {\cal H}^r$, $r > s > 0$, and $f_\eps \in H$ be such that 
$\|f-f_\eps\|_H^2 \leq \eps$, where $\eps > 0$. Then there exists a smoothing 
operator $S: H \to {\cal H}^s$ and a positive function $\gamma$ with 
$\lim_{x\downarrow 0} \gamma(x) = 0$, such that $\tilde{f}_\eps := 
S f_\eps \in {\cal H}^s$ satisfies $\|f-\tilde{f}_\eps\|_s^2 \leq 
\gamma(\eps)$. \\
Proof:
Using the resolution of the identity associated to $A$, we define for $h>0$ 
the operator $S_h: H \to {\cal H}^s$ by $S_h := \int_0^{1/h} d E_\lbd$.%
\footnote{Recall that $\int_0^\infty d E_\lbd$ is the identity operator 
in $H$.}
Defining $\tilde{f}_\eps := S_h f_\eps \in {\cal H}^s$ one can estimate
\begin{eqnarray}
\| f - \tilde{f}_\eps \|_s^2
   &  =   & \| f - S_h f_\eps \pm S_h f \|_s^2 \nonumber \\
   & \leq & 2 \left( \| (I-S_h)f \|_s^2 + \| S_h(f-f_\eps) \|_s^2 \right) . 
              \label{glaet-abschaetz1}
\end{eqnarray}
The first term on the right hand side of (\ref{glaet-abschaetz1}) can be 
estimated by
\begin{eqnarray*}
\| (I-S_h)f \|_s^2
   &  =   & \int_{1/h}^\infty (1+\lbd^2)^s \left[ \frac{(1+\lbd^2)^r}
            {(1+\lbd^2)^r} \right] d \ipl E_\lbd f, f \ipr \\
   & \leq & (1+h^{-2})^{s-r} \int_0^\infty (1+\lbd^2)^r
            d \ipl E_\lbd f, f \ipr \\
   &  =   & (1+h^{-2})^{s-r} \| f \|_r^2 .
\end{eqnarray*}
For the second term on the right hand side of (\ref{glaet-abschaetz1}) we have
\begin{eqnarray*}
\| S_h(f-f_\eps) \|_s^2
   &  =   & \int_0^{1/h} (1+\lbd^2)^s d\ipl E_\lbd (f-f_\eps),(f-f_\eps)\ipr \\
   & \leq & (1+h^{-2})^s \int_0^{1/h} d\ipl E_\lbd (f-f_\eps),(f-f_\eps)\ipr \\
   & \leq & (1+h^{-2})^s \| f-f_\eps \|_H^2 .
\end{eqnarray*}
Substituting the last inequalities in (\ref{glaet-abschaetz1}) we obtain
\begin{equation} \label{glaet-abschaetz2}
\| f - \tilde{f}_\eps \|_s^2  \leq  
   2 \left[ (1+h^{-2})^s \eps + (1+h^{-2})^{s-r} \| f \|_r^2 \right] .
\end{equation}
To balance the right hand side of (\ref{glaet-abschaetz2}) one must choose 
$h = [ (\eps^{-1}\| f \|_r^2)^{1/r} - 1]^{-\me}$. Now the theorem follows 
choosing $S := S_h$ and $\gamma(x) := 4 x^{(r-s)/r}\| f \|_r^{2s/r}$. \eop
\end{lemma}

\begin{remrk}
        Let $f_\eps \in H$ be the given Cauchy data for $(P_e)$. From 
Lemma~\ref{lemma-fehler-glaetung} follows that when the exact Cauchy data 
$f$ is better than ${\cal H}^{\me}$, i.e. $f \in {\cal H}^r$ for 
$r > \me$, then it is possible to find a $\tilde{f}_\eps$ in ${\cal H}^{\me}$ 
near to $f$ in the $(\me)$--norm. For the hyperbolic case one obtains an 
analogous result. \eop
\end{remrk}

        Since the affine term $z_{f,g}$ depends continuously on the data 
$(f,g)$ -- respectively $z_f$ depends continuously on $f$, we conclude 
from Lemma~\ref{lemma-fehler-glaetung} that under the corresponding 
assumption it is possible to obtain from the measured data $(f_\eps,g_\eps) 
\in H^2$ a $z_\eps$ satisfying $\|z_{f,g} - z_\eps\| \leq \eps'$ 
(respectively $\|z_f - z_\eps\| \leq \eps'$).

        Let $T$ be one of the operators defined in (\ref{Te-oper-def}), 
(\ref{Th-oper-def}) or (\ref{Tp-oper-def}) and $T_l$ the corresponding linear 
part. We want to choose a linear operator $R$ such that given $\vphi_0$, the 
regularized sequence $\tilde{\vphi}_{k+1} := R \tilde{\vphi}_k + z_\eps$ 
converges faster then the original one $\vphi_{k+1} := T_l \vphi_k + z_{f,g}$. 
Simultaneously we have to assure that the difference $\|\lim \tilde{\vphi}_k 
- \lim \vphi_k\|$ remains small.

        In section~\ref{sec-methoden} we have seen that $T_l = \int_0^\infty 
F(\lbd) d E_\lbd$, where $F(\lbd)$ is either $\tanh^2(\lbd)$, $\cos^2(\lbd)$ 
or $(1 - \gamma \exp(-\lbd^2T))$. Given $n \in \N$ we define the 
regularization operator $R_n$ by
\[ R_n := \int_0^n F(\lbd) d E_\lbd . \]
Next we define $\bar{\vphi}$ and $\vphi_n$ as the fixed points of $\bar{\vphi} 
= T_l \bar{\vphi} + z_{f,g}$ and $\vphi_n = R \vphi_n + z_\eps$ respectively 
(note that $\vphi_n$ exists since $R_n$ is contractive). From the identity
\[ \vphi_n - \bar{\vphi} = R_n (\vphi_n - \bar{\vphi}) +
   (R_n - T_l) \bar{\vphi} + z_\eps - z_{f,g} \]
one obtains the estimate
\begin{equation} \label{dirkrep-gl}
\| \vphi_n - \bar{\vphi} \| \leq \| (I - R_n)^{-1}(R_n - T_l) \bar{\vphi} \|
   + \eps' \| (I - R_n)^{-1} \| ,
\end{equation}
which leads us to the following lemma

\begin{lemma} \label{diskrep-lemma}
        Let $T_l$ represent the linear part of the iterative procedure for one 
of the problems $(P_e)$, $(P_h)$ or $(P_p)$. Given the corresponding family of 
operators $R_n$ defined as above, we have
\[ \lim_{n\to\infty} \| (I - R_n)^{-1}(R_n - T_l) \bar{\vphi} \| = 0 \ \ \ \
{\rm and} \ \ \ \ \lim_{n\to\infty} \| (I - R_n)^{-1} \| = \infty . \]
Proof:
Since $(I - R_n)$ is the identity operator on Rg$(R_n - T_l)$, the first 
assertion follows from the inequality%
\footnote{Here $s \in \R$ must be chosen according to the space where the 
iteration takes place.}
\[ \| (I - R_n)^{-1}(R_n - T_l) \bar{\vphi} \|_s^2  \leq
   \int_n^\infty (1+\lbd^2)^s d \ipl E_\lbd \bar{\vphi} , \bar{\vphi} \ipr . \]
The second assertion follows from the identity $\|(I - R_n)^{-1}\| = 
(1-\Lambda(n))^{-1}$, where $\Lambda(n):= \max\{ \lbd\in\sigma(A); \lbd<n \}$, 
and the fact that $A: H \to H$ is unbounded. \eop
\end{lemma}

        Now making {\em a priori} assumptions on the regularity of 
$\bar{\vphi}$, we obtain from Lemma~\ref{diskrep-lemma} and the estimate 
(\ref{dirkrep-gl}) the desired regularization result.

\begin{lemma} \label{regulariz-lemma}
        If there exists a positive monotone increasing function $G \in 
C(\R^+)$ with 
\[ \lim_{\lbd\to\infty} G(\lbd) = \infty \ \ \ \ {\rm and} \ \ \ \
   \int_0^\infty (1+\lbd^2)^s G^2(\lbd) d \ipl E_\lbd 
   \bar{\vphi} , \bar{\vphi} \ipr = M^2 < \infty , \]
then exists an optimal choice of $n^* \in \N$ such that $\|\vphi_{n^*} - 
\bar{\vphi}\| \leq \|\vphi_n - \bar{\vphi}\|$ for all $n \in \N$. Further 
$n^*$ solves the minimization problem
\[ \min_{n\in\N} \left\{ M G^{-1}(n) + \eps (1+\Lambda(n))^{-1} \right\} . \]
Proof:
From Lemma~\ref{diskrep-lemma} follows
\begin{eqnarray*}
\| (I - R_n)^{-1}(R_n - T_l) \bar{\vphi} \|_s^2  & \leq &
   \int_n^\infty (1+\lbd^2)^s \frac{G^2(\lbd)}{G^2(\lbd)}
   d \ipl E_\lbd \bar{\vphi} , \bar{\vphi} \ipr \\
 & \leq & G^{-2}(n) \int_n^\infty (1+\lbd^2)^s G^2(\lbd)
   d \ipl E_\lbd \bar{\vphi} , \bar{\vphi} \ipr \\
 & \leq & M^2 G^{-2}(n) .
\end{eqnarray*}
From (\ref{dirkrep-gl}) we obtain
\[ \| \vphi_n - \bar{\vphi} \| \leq M G^{-1}(n) + \eps (1+\Lambda(n))^{-1} \]
and the theorem follows. \eop
\end{lemma}
%
%
%
% --------------------------------------------------------------------------- %
%                                   Section 6                                 %
% --------------------------------------------------------------------------- %
\setcounter{footnote}{1}
\section{Numerical results} \label{sec-numerik}
%
%
%
%----------------------------------------------------------------------------%
\subsection{A parabolic reconstruction problem}

        We consider the heat equation $a^2 \partial_t u = \Delta u$ at 
$(0,T)\times\Omega$, where $\Omega = (0,1)\times(0,1)$. The solution $u(0)$ 
of the reconstruction problem is shown in Figure~1. It consists of a 
$L^2(\Omega)$ function added to a polynomial of fourth degree.

\begin{figure}                                                        %figure 1
\centerline{ \epsfysize5cm \epsfxsize7.5cm \epsfbox{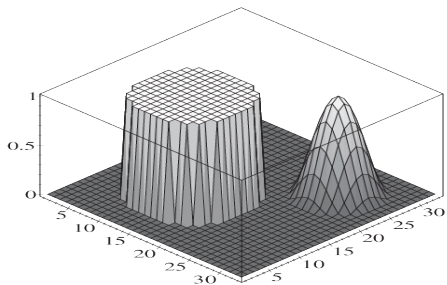} }
\vskip-0.5cm
\caption{$u(0)$: solution of the reconstruction problem}
\end{figure}

        In the first example we choose $a^2 = 8$ and take as problem data 
$f := u(T)$ evaluated at $T = 0.0625$. The iterative procedure is started 
with $\vphi_0 \equiv 0$ and we chose the parameter $\gamma=2$, which is in 
agreement with Lemma~\ref{lemma-Tlp-eigenschaft}. In Figure~2 one can see 
the data $f$ of the reconstruction problem and the iteration error after 
$10^6$ steps.

                                                                      %figure 2
\centerline{ \epsfysize5cm \epsfxsize7.5cm \epsfbox{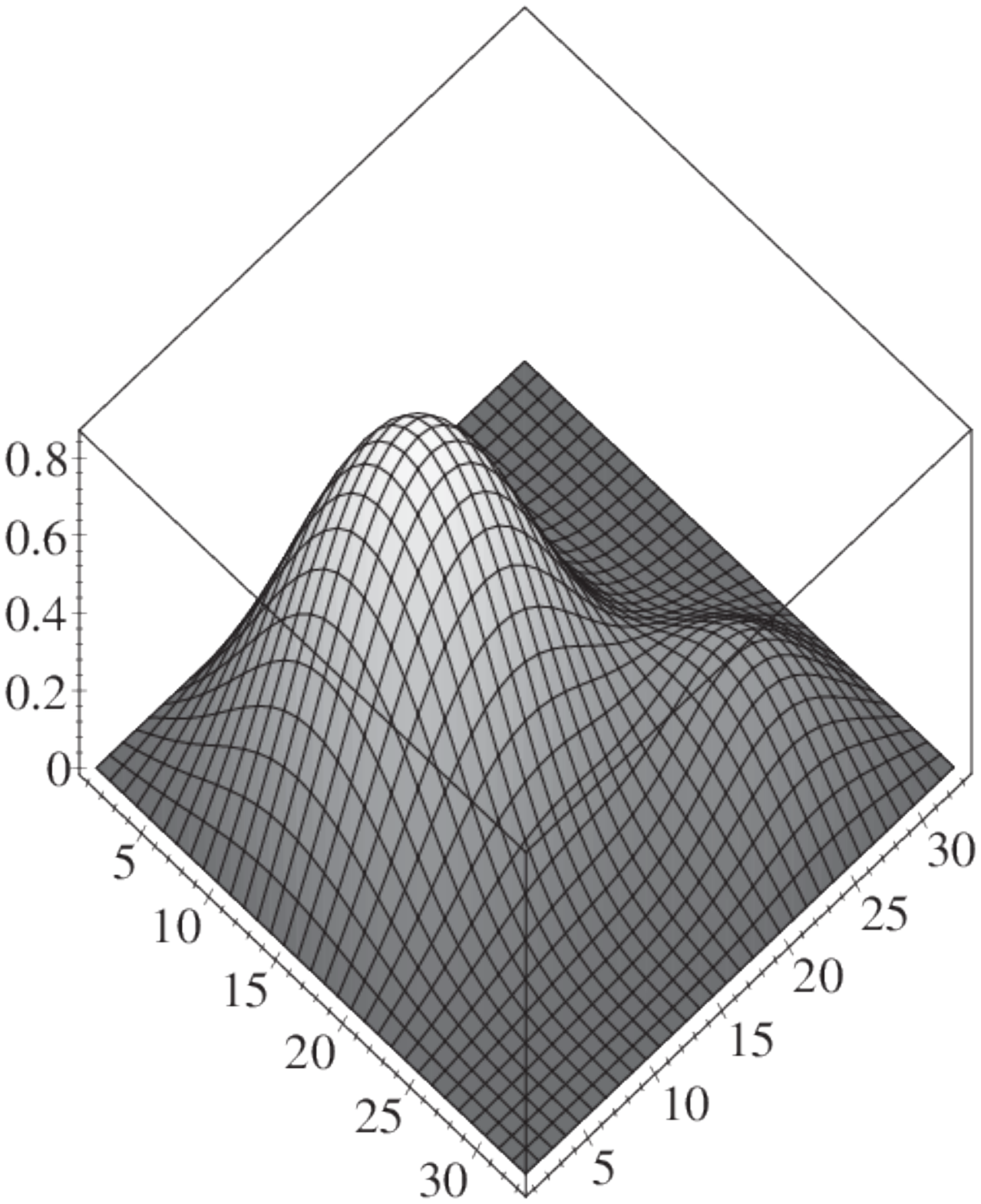}
             \epsfysize5cm \epsfxsize7.5cm \epsfbox{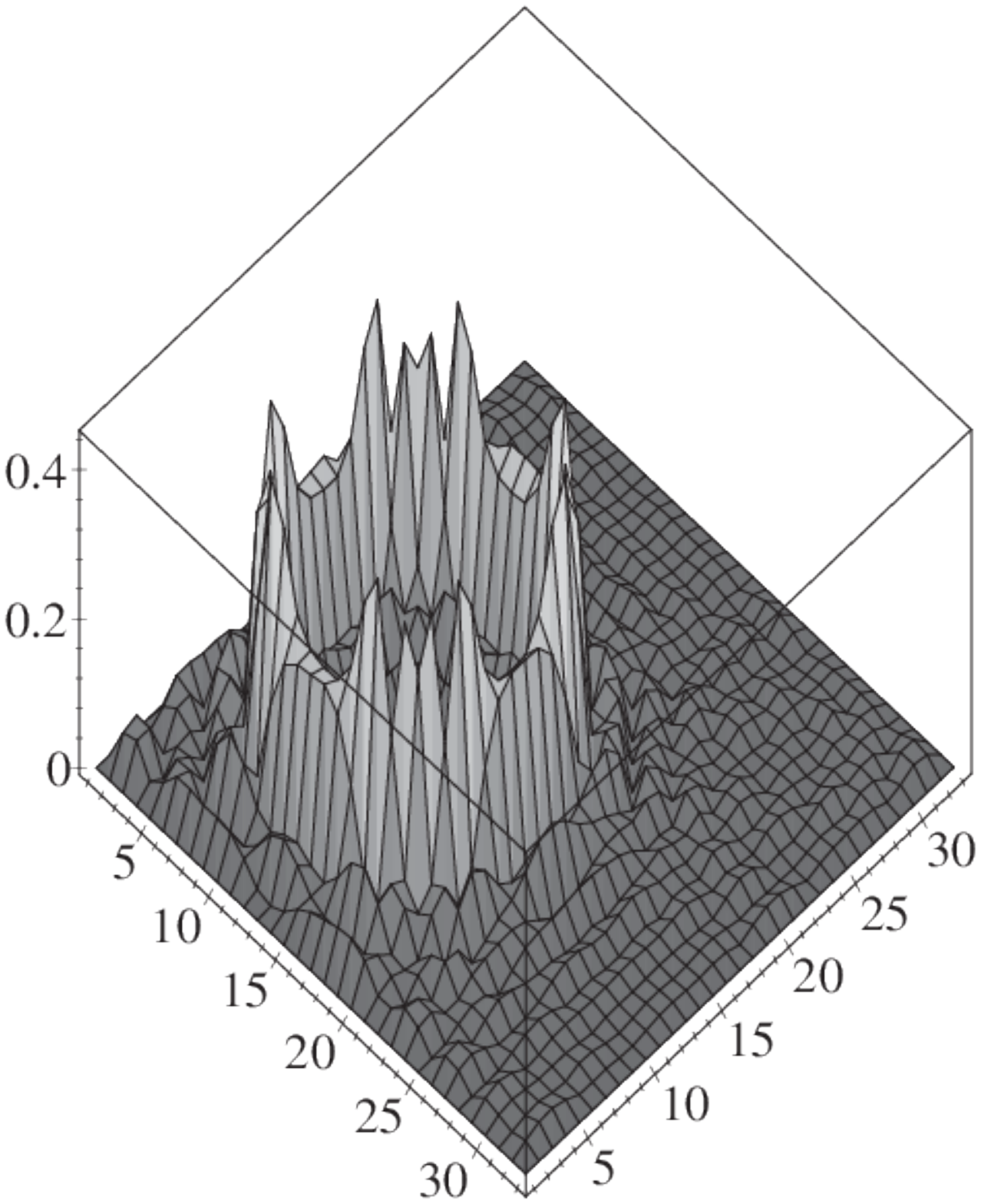}  }
\vskip-0.2cm
\centerline{Figure~2: Temperature profile $u(T)$ for $a^2=8$ and iteration 
error after $10^6$ steps}

        In the second example we choose $a^2 = 2$ and set $f := u(T)$ for 
$T = 0.0625$, where $u(0)$ is the same as before. The iterative procedure is 
started with $\vphi_0 \equiv 0$. In Figure~3 one can see the problem data 
$f$ and the iteration error after $10^6$ steps.

\setcounter{figure}{2}
\begin{figure}                                                        %figure 3
\centerline{ \epsfysize5cm \epsfxsize7.5cm \epsfbox{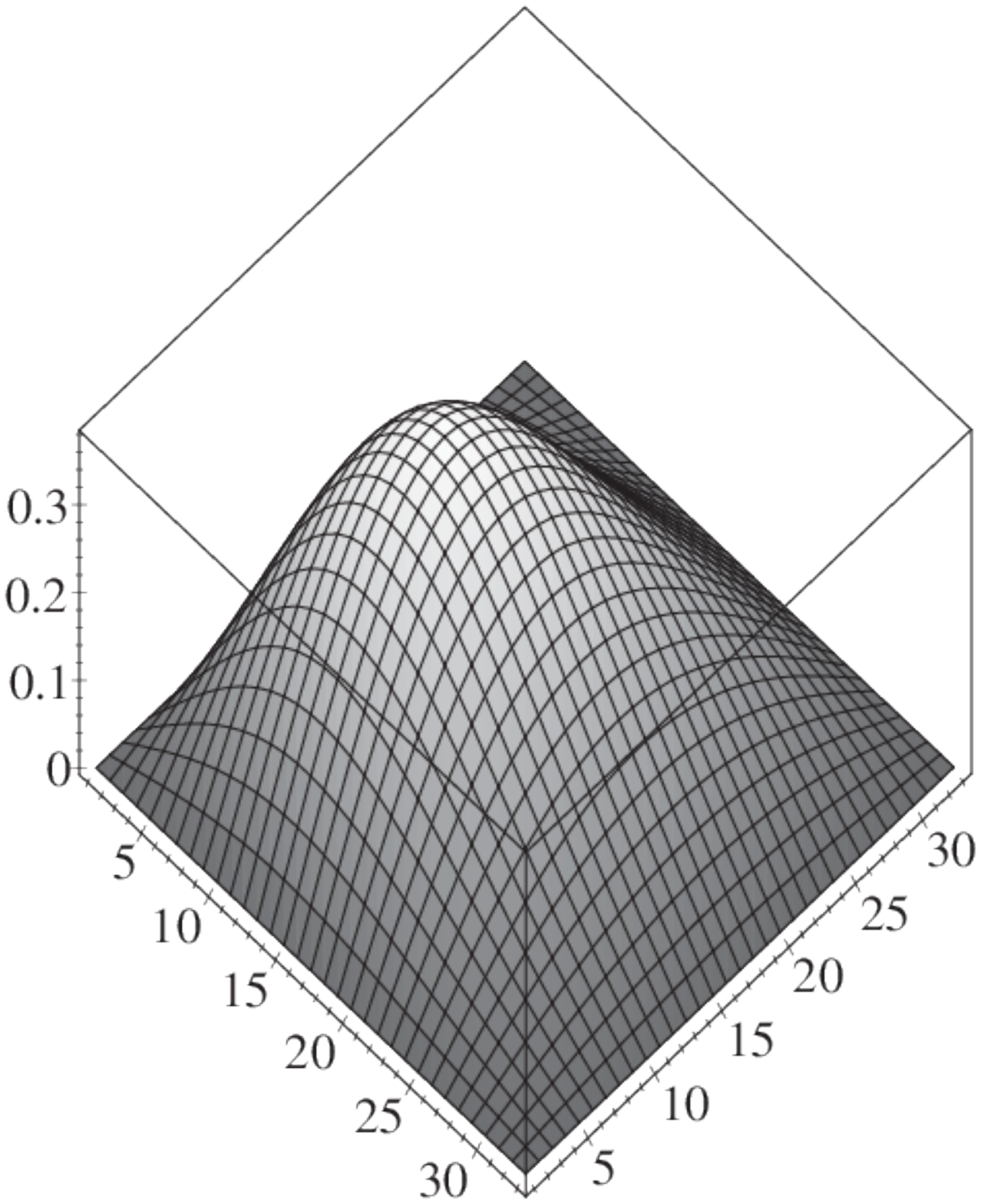}
             \epsfysize5cm \epsfxsize7.5cm \epsfbox{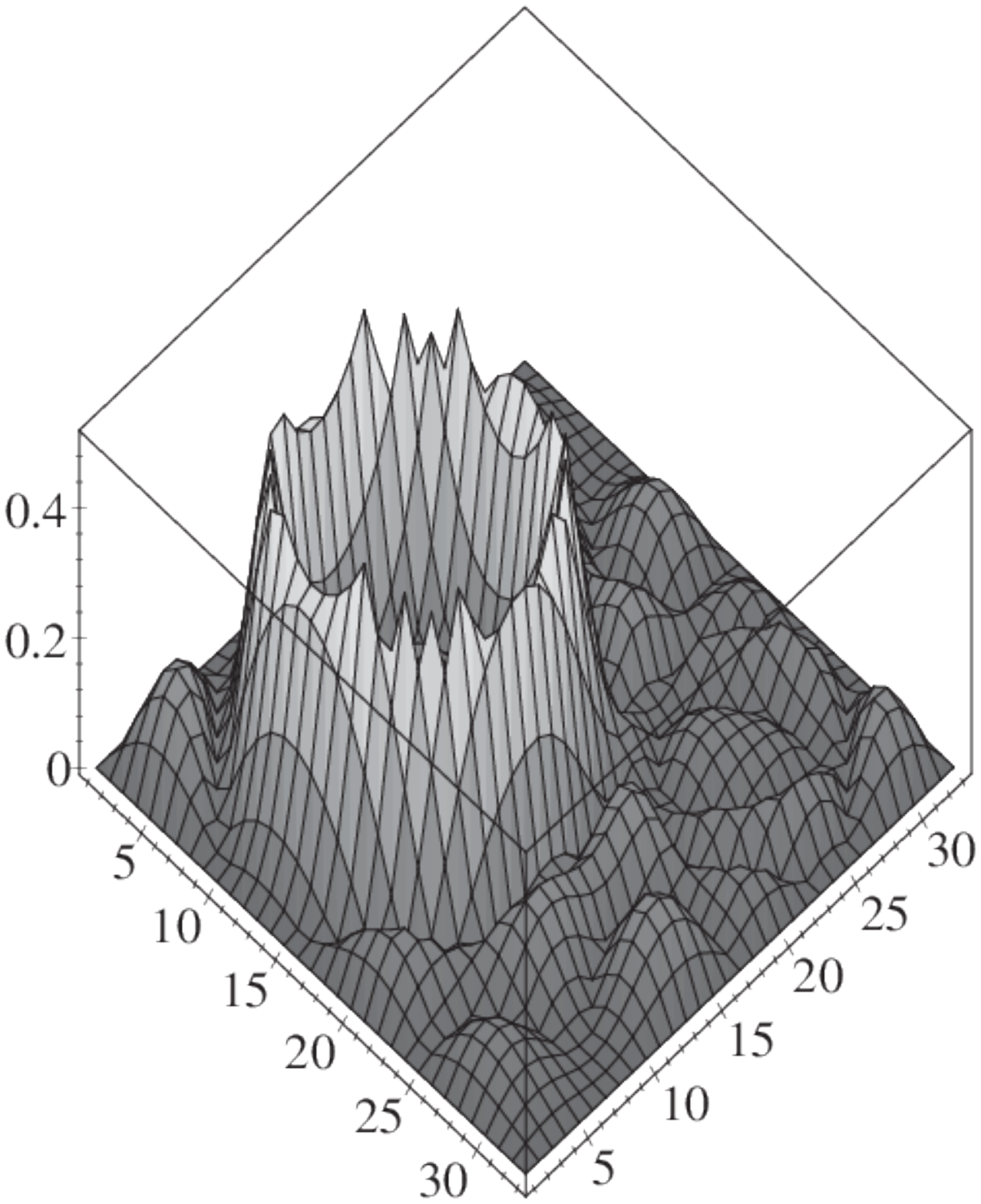} }
\vskip-0.5cm
\caption{Temperature profile $u(T)$ for $a^2=2$ and iteration error after 
$10^6$ steps}
\end{figure}

        One should note that $\bar{\vphi} := u(0)$ is a fixed point of 
the numerical iteration. This follows from the fact that $f  = u(T)$ was 
obtained by solving a direct problem.

        In both examples the reconstruction is much better at the part of 
the domain where the initial condition is smooth. In Table~1 we present the 
evolution of the iteration error $\vphi_k - u(0)$ for the two examples above.

        Note also that the convergence speed decays exponentially as we 
iterate. This is a consequence of the exponential behavior of the eigenvalues 
of $T_{l,p}$ (see Paragraph~\ref{ssec-parab-analyse}). \\[2ex]
\centerline{Table~1: Evolution of the relative error in the $L^2$--norm 
(parabolic problem)} \\[1ex]
\centerline{
\begin{tabular}{cccccc} \hline
{} & $10$ steps & $10^3$ steps & $10^4$ steps & $10^5$ steps & $10^6$ steps \\
\hline
$a^2 = 8$ & 32.5\%  & 25.7\% & 23.5\% & 22.5\% & 20.9\%  \\ \hline
$a^2 = 2$ & 49.8\%  & 42.2\% & 40.1\% & 36.2\% & 31.4\%  \\ \hline
\end{tabular} }
%
%
%
%----------------------------------------------------------------------------%
\subsection{An elliptic reconstruction problem}

        We consider next the Laplace equation $\partial^2_t u + \partial^2_x u 
= 0$ at $(0,T)\times\Omega$, where $\Omega=(0,1)$. Given $k \in \N$ we choose 
the Cauchy data $f \equiv 0$, $g_k = \sin(k\pi x)$ and try to reconstruct the 
corresponding traces $\partial_t u(T)$ at the final time $T=1$.%
\footnote{Note that $g_k = \sin(k\pi x)$ are eigenfunctions of $T_{l,e}$ with 
corresponding eigenvalues $\lbd_k = \tanh(k\pi)^2$. The solutions of the 
reconstruction problems are given by $\cosh(k\pi t) \sin(k\pi x)$.} \\[2ex]
\centerline{Table~2: Evolution of the relative error in the $L^2$--norm 
(elliptic problem)} \\[1ex]
\centerline{
\begin{tabular}{ccccccc} \hline
{} & $10^2$ steps & $10^3$ steps & $10^5$ steps
   & $10^6$ steps & $10^8$ steps & $10^9$ steps \\
\hline
$k = 1$ &    48.33\% & \,\ 5.34\% & \,\ 5.30\% & \,\ 5.30\% & \,\ 5.30\%
        & \,\ 5.30\% \\ \hline
$k = 2$ &    99.86\% &    98.61\% &    29.72\% &    11.28\% &    11.28\%
        &    11.28\% \\ \hline
$k = 3$ &    99.99\% &    99.98\% &    99.74\% &    97.44\% &    21.79\%
        &    18.42\% \\ \hline
\end{tabular} } \\[2ex]
        In Table~2 the evolution of the relative reconstruction error for 
three distinct values of $k$ is presented. From this data one can see that if 
$g$ can be expanded in a Fourier series, it's first coefficient will be 
accurately reconstructed after $10^4$ steps, while the second one only after 
$10^6$ steps, etc \dots

\section*{Acknowledgments}
Research partially supported by CNPq/GMD under grant 91.0206/98-8.


\begin{thebibliography}{99999}
\addcontentsline{toc}{chapter}{Literaturverzeichnis}

\bibitem[Bas]{Bas} {\sc Bastay, G.,} {\it Iterative Methods for Ill--Posed
                Boundary value Problems,} Link\"oping Studies in Science and
                Technology, Dissertations No. {\bf 392}, Link\"oping, 1995.

\bibitem[Bau]{Bau} {\sc Baumeister, J.,} {\it Stable Solution of Inverse
                Problems,} Fried.Vieweg \& Sohn, Braunschweig, 1987.

\bibitem[Is]{Is} {\sc Isakov, V.M.,} {\it On the Uniqueness of the Solution of
                the Cauchy Problem,} Soviet Math. Dokl., Vol {\bf 22} (1980),
                No. {\bf 3}, 639--642.

\bibitem[Je]{Je} {\sc Jeggle, H.,} {\it Nichtlineare Funktionalanalysis,}
                Teubner, Stuttgart, 1979.

\bibitem[JoNa]{JoNa} {\sc Jourhmane, M. and Nachaoui, A.,} {\it A Relaxation
                Algorithm for Solving a Cauchy--Problem,} Preliminary
                Proceedings--Vol {\bf 2}, 2nd Intern. Confer. on Inverse
                Problems in Engineering: Theory and Practice, Le Croisic,
                June 1996.

\bibitem[KM1]{KM1} {\sc Kozlov, V.A., Maz'ya, V.G. and Fomin, A.V.,} {\it An
                iterative method for solving the Cauchy problem for elliptic
                equations,} Comput.Maths.Phys., Vol. {\bf 31} (1991), No. 
                {\bf 1}, 45--52.

\bibitem[KM2]{KM2} {\sc Kozlov, V.A. and Maz'ya, V.G.,} {\it On iterative 
                procedures for solving ill-posed boundary value problems that
                preserve differential equations,} Leningrad Math. J., Vol.
                {\bf 1} (1990), No. {\bf 5}, 1207--1228.

\bibitem[Kra]{Kra} {\sc Krasnosel'skii, M.A., Vainikko, G.M., Zabreiko, P.P.,
                Rutitskii, Yu B. and Stetsenko, V.Yu.,} {\it Approximate
                Solution of Operator Equations,} Wolters--Nordhoff Publishing,
                Groningen, 1972.

\bibitem[Le1]{Le1} {\sc Leit\~ao, A.,} {\it Ein Iterationsverfahren f\"ur
                elliptische Cauchy--Probleme und die Ver\-kn\"upfung mit der
                Backus--Gilbert Methode,} Dissertation, FB Mathematik, J.W. 
                Goethe--Universit\"at, Frankfurt am Main, 1996.

\bibitem[Le2]{Le2} {\sc Leit\~ao, A.,} {\it An Iterative Method for Solving
                Elliptic Cauchy Problems,} Numerical Functional Analysis and
                Optimization, to appear.

\bibitem[LiMa]{LiMa} {\sc Lions, J.L. and Magenes, E.,} {\it Non--Homogeneous
                Boundary Value Problems and Applications,} Springer--Verlag,
                Berlin Heidelberg New York, 1972.

\bibitem[Va]{Va} {\sc Vainikko, G.M.,} {\it Regularisierung nichkorrekter
                Aufgaben,} Preprint No. {\bf 200}, Universit\"at
                Kaiserslautern, 1991.

\end{thebibliography}
\end{document}